# A 'replacement sequence' method for finding the largest real root of an integer monic polynomial

A. K. Gupta, Department of Electronics and Communication, Allahabad University, Allahabad - 211 002, India
(e-mail: ashok@mri.ernet.in)

A. K. Mittal, Department of Physics, Allahabad University, Allahabad – 211 002, India
(e-mail: mittal@mri.ernet.in)

ABSTRACT:

To every integer monic polynomial of degree *m* can be associated a 'replacement rule' that generates a word **W*** from another word **W** consisting of symbols belonging to a finite 'alphabet' of size 2*m*. This rule applied iteratively on almost any initial word $\mathbf{W_0}$, yields a sequence of words $\{\mathbf{W_i}\}$. From a count of different symbols in the word $\mathbf{W_i}$, one can obtain a rational approximate to the largest real root of the polynomial.

Let $p(x) = x^m - a_1 x^{m-1} - \ldots - a_m$ , where $a_i$ are integers, be the given monic polynomial.

Let $A^+ = \{A^+_1, A^+_2, \ldots, A^+_m\}$ be a finite set of *m* symbols. Let $A^\sim = \{A^\sim_1, A^\sim_2, \ldots, A^\sim_m\}$ be another set of *m* symbols associated with the set A so that with each element $A^+_i$ of $A^+$ is associated the element $A^\sim_i$ of $A^\sim$. Let $A = A^+ \cup A^\sim$. Elements of A are called 'letters' belonging to the 'alphabet' A. Let $A^p = A \times A \times \ldots \times A$ (p times). Elements of $A^p$ are called 'words' of length p that can be formed from 'letters' of A. Let $A^* = \cup_p A^p$. Clearly, $A^*$ consists of all 'words' that can be made from the 'alphabet' A. A 'letter' $A^+_i$ repeated k times in a word will be denoted by $(A^+_i)^k$. Similarly a 'letter' $A^\sim_i$ repeated k times in a word will be denoted by $(A^\sim_i)^k$. We define $(\alpha^+_i)^k$ by

$$
\begin{aligned}
(\alpha^+_i)^k &= (A^+_i)^k && \text{if } k \geq 0 \\
(\alpha^+_i)^k &= (A^\sim_i)^k && \text{if } k < 0 \qquad (1)\\
(\alpha^\sim_i)^k &= (A^\sim_i)^k && \text{if } k \geq 0 \\
(\alpha^\sim_i)^k &= (A^+_i)^k && \text{if } k < 0
\end{aligned}
$$

Consider the replacement rule R: $A \rightarrow A^*$ given by

$$
\begin{aligned}
A^+_1 &\rightarrow (\alpha^+_1)^{a1} A^+_1 A^+_2 \\
A^+_2 &\rightarrow (\alpha^+_1)^{a2} A^+_2 A^+_3 \\
&. \\
A^+_i &\rightarrow (\alpha^+_1)^{ai} A^+_i A^+_{i+1} \\
&. \\
&. \\
A^+_m &\rightarrow (\alpha^+_1)^{am} A_m \\
A^~_1 &\rightarrow (\alpha^~_1)^{a1} A^~_1 A^~_2 \\
A^~_2 &\rightarrow (\alpha^~_1)^{a2} A^~_2 A^~_3 \\
&. \\
A^~_i &\rightarrow (\alpha^~_1)^{ai} A^~_i A^~_{i+1} \\
&. \\
&. \\
A^~_m &\rightarrow (\alpha^~_1)^{am} A_m
\end{aligned}
\qquad (2a)
$$

Let $W = A_{s1}A_{s2}\ldots\ldots A_{sq}$, where sk (k=1,2,……..,q) $\in$ {1,2,……..,m}. Then $W \in A^*$. The replacement rule R induces a mapping $R^*:A^* \rightarrow A^*$ defined by

$$
\begin{aligned}
W^* &= R^*(W) \\
&= R^*(A_{s1}A_{s2}\ldots\ldots A_{sq}) \\
&= R(A_{s1})R(A_{s2})\ldots\ldots R(A_{sq})
\end{aligned}
\qquad (2b)
$$

Let $W_0 \in A^*$. Then

$$W_i = (R^*)^i (W_0) = R^*((R^*)^{i-1}(W_0)) = R^*(W_{i-1}) \qquad (3)$$

denotes the word obtained by i times repeated application of the replacement rule (2) on the initial word $W_0$.

Let $\mathbf{n}: A^* \rightarrow V^m$ be a mapping which assigns to a word in A* a vector

$$\mathbf{n}(W) = (n^+_1(W) - n^~_1(W), n^+_2(W) - n^~_2(W), \ldots\ldots, n^+_i(W) - n^~_i(W), \ldots\ n^+_m(W) - n^~_m(W))^T \qquad (4)$$

in an m-dimensional vector space $V^m$ such that $n^+_i(W)$ and $n^~_I(W)$ are non-negative integers denoting the frequency of occurrence of $A^+_j$ and $A^~_i$ in the word W.

The replacement rule (2) and the mapping **n** induces a (m x m) matrix **R** which maps vectors in $V^m$ into $V^m$ such that

$$\mathbf{n}(W^*) = \mathbf{n}(R^*(W)) = \mathbf{R}(\mathbf{n}(W)) \tag{5}$$

It follows from (2) - (5) that

$$\mathbf{R} = \begin{bmatrix} 1+a_1 & a_2 & a_3 & \ldots & a_{n-1} & a_n \\ 1 & 1 & 0 & \ldots & 0 & 0 \\ 0 & 1 & 1 & \ldots & 0 & 0 \\ . & . & . & \ldots & . & . \\ . & . & . & \ldots & 1 & . \\ 0 & 0 & 0 & \ldots & 1 & 1 \end{bmatrix} \tag{6}$$

Equations (3) and (5) imply

$$\mathbf{n}(W_i) = \mathbf{n}(R^*(W_{i-1})) = \mathbf{R}(\mathbf{n}(W_{i-1})) = \mathbf{R}^i(\mathbf{n}(W_0)) \tag{7}$$

For almost any $W_0$, $\mathbf{n}(W_i)$ tends, as $i \to \infty$, to multiples of the eigenvector corresponding to the eigenvalue with maximum absolute value [1]. The matrix R in equation (6) is equal to the identity matrix plus the companion matrix of the polynomial p. The eigenvalues of the matrix **R** are 1 plus the roots of the polynomial p . The eigenvector corresponding to the largest real eigenvalue is given by $[\lambda^{m-1}, \lambda^{m-2}, \ldots, \lambda, 1]^T$, where $\lambda$ is the largest real root of the polynomial p. Thus one finds that for almost any initial word $W_0$,

$$\lim_{i \to \infty} \quad n_j(W_i)/n_{j+1}(W_i) = \lambda \qquad j = 1,2, \ldots, m-1 \tag{8}$$

Eqn (8) shows that the largest real root of any integer monic polynomial can be obtained to any accuracy only by operations consisting of replacement of symbols, counting and subtracting. This method does not require any ‘advanced’ operation such as multiplication, division, logarithms and anti-logarithms, although the justification of the method uses ‘advanced’ concepts like the theory of matrices.

The numbers that can be obtained geometrically using only a ruler and a compass are called constructible numbers [2]. It is known that the set of constructible numbers is $Q(\sqrt{m_1}, \sqrt{m_2}, \sqrt{m_3}, \ldots, \sqrt{m_N})$, the extension field of rational numbers by square roots of integers. Higher roots such as cube roots etc, in general, are not constructible. The replacement and counting method of our earlier paper[3] shows how to construct these roots by a method more primitive than the use of ruler and compass. The set of numbers that may be constructed by this method, which may be called *constructible by replacement and counting,* contains the set of (geometrically) constructible numbers $Q(\sqrt{m_1}, \sqrt{m_2}, \ldots, \sqrt{m_N})$. This paper generalizes the earlier paper[3] to show that the largest real root of a integer monic polynomial can be constructed by *replacement, counting and subtracting.*